\documentclass[10pt,openany,leqno]{amsart}
\usepackage[backref=page,colorlinks=true]{hyperref}
\usepackage[capitalize]{cleveref} 
\usepackage[english, francais]{babel}
\usepackage{amsmath,amsthm,amsfonts,amssymb,amscd,url}
\usepackage{graphics}
\numberwithin{equation}{section}
\input{xy} \xyoption{all}
\usepackage[latin1]{inputenc}

\usepackage{enumerate, color}
 
\usepackage{hyperref, mathrsfs}
\usepackage{epsfig} 
\input{xy} \xyoption{all}

\usepackage[scr=boondoxo,scrscaled=1.05]{mathalfa}

\begin{document}
\def é{\' e}
\def \cE{\mathscr E}
\def\inv{^{-1}}
\def \Card{\mathrm{Card\, }}
\def \Leb{\mathrm{Leb\, }}
\def\Max{\text{max\, }}
\def\cal{\mathcal}
\def\R{\mathbb R}
\def\N{\mathbb N}
\def\Z{\mathbb Z}
\def\C{\mathbb C}
\def\cC{\mathcal C}
\def\T{\mathbb T}
\def\a{{\underline a}}
\def\b{{\underline b}}
\def\c{{\underline c}}
\def\Log{\text{log}}
\def\loc{\text{loc}}
\def\inta{\text{int }}
\def\det{\mathrm{det\, }}
\def\exp{\mathrm{exp\, }}
\def\Re{\text{Re}}
\def\lip{\mathrm{Lip\, }}
\def\leb{\mathrm{Leb\, }}
\def\dom{\text{Dom}}
\def\diam{\mathrm{diam}\:}
\def\supp{\mathrm{supp}\:}
\newcommand{\ovfork}{{\overline{\pitchfork}}}
\newcommand{\ovforki}{{\overline{\pitchfork}_{I}}}
\newcommand{\Tfork}{{\cap\!\!\!\!^\mathrm{T}}}
\newcommand{\whforki}{{\widehat{\pitchfork}_{I}}}
\newcommand{\marginal}[1]{\marginpar{{\scriptsize {#1}}}}
\def \np{{\color{red} 4}}
\def \npm{{\color{red} 3}}
\def \npmm{{\color{red} 2}}
\def\cP{{\mathscr P}}
\def\cF{{\mathscr F}}
\def\cH{{\mathcal H}}
\def\cV{{\mathcal V}}
\def\cX{{\mathcal X}}
\def\cY{{\mathcal Y}}
\def\qand{\quad \text{and} \quad}
\def\sR{{\mathfrak R}}
\def\sS{{\mathfrak S}}
\def\sU{{\mathfrak U}}
\def\sM{{\mathfrak M}}
\def\sA{{\mathfrak A}}
\def\sB{{\mathfrak B}}
\def\sX{{\mathfrak X}}
\def\sY{{\mathfrak Y}}
\def\sZ{{\mathfrak Z}}
\def\sE{{\mathfrak E}}
\def\sP{{\mathfrak P}}
\def\sG{{\mathfrak G}}
\def\sa{{\mathfrak a}}
\def\sb{{\mathfrak b}}
\def\sc{{\mathfrak c}}
\def\sd{{\mathfrak d}}
\def\se{{\mathfrak e}}
\def\sg{{\mathfrak g}}
\def\sh{{\mathfrak h}}
\def\sn{{\mathfrak n}}
\def\sq{{\mathfrak q}}
\def\sp{{\mathfrak p}}
\def\sl{{\mathfrak {l}}}
\def\Diff{\mathrm {Diff}}
\def\sr{{\mathfrak {r}}}
\def\ss{{\mathfrak {s}}}
\def\st{{\mathfrak {t}}}
\def\su{{\mathfrak {u}}}
\def\spp{{\mathfrak {p}}}
\def\sw{{\mathfrak {w}}}
\def\arr{\overleftarrow}
\def\avv{\overrightarrow}
\theoremstyle{plain}
\newtheorem{theo}{\bf Theorem}[section]
\newtheorem{theorem}[theo]{\bf Theorem}
\newtheorem{lemma}[theo]{\bf Lemma}
\newtheorem{conjecture}[theo]{\bf Conjecture}
\newtheorem{prop}[theo]{\bf Proposition}
\newtheorem{problem}[theo]{\bf Problem}
\newtheorem{question}[theo]{\bf Question}
\newtheorem{proposition}[theo]{\bf Proposition}
\newtheorem{claim}[theo]{\bf Claim}
\newtheorem{assumption}[theo]{\bf Assumption}
\newtheorem{coro}[theo]{\bf Corollary}
\newtheorem{corollary}[theo]{\bf Corollary}

\newtheorem{sublemma}[theo]{\bf Sublemma}
\newtheorem{rema}[theo]{\bf Remark}
\newtheorem{remark}[theo]{\bf Remark}
\newtheorem{fact}[theo]{\bf Fact}

\newtheorem{exem}[theo]{\bf Example}
\newtheorem{definition}[theo]{\bf Definition}
\newtheorem{defi}[theo]{\bf Definition}


\title{
Complexities of differentiable dynamical systems\\
}

\author{Pierre Berger}
\address{Institut de Mathématiques de Jussieu-Paris Rive Gauche\\
CNRS, UMR 7586\\
Sorbonne université\\
4 place Jussieu, 75252 Paris Cedex 05}
\maketitle
\begin{otherlanguage}{english}
\begin{abstract}
We define the notion of localizable property for a dynamical system. 
Then we survey three properties of complexity and relate how they are known to be typical among differentiable dynamical systems. These notions are the fast growth of the number of periodic points, the positive entropy and the high emergence. 
We finally propose a dictionary between the previously explained theory on entropy and the ongoing one on emergence.
\end{abstract}
\section*{Introduction}

A differentiable dynamical systems is a  map $f\in C^r(M,M)$ of a manifold $M$, for $1\le r\le \infty$. Our aim is to understand the properties of the  orbits $(f^n(x))_{n\ge 0}$ of the points $x\in M$, where $f^n(x)=f\circ \cdots \circ f(x)$. Let us start from the abstract viewpoint. A property $(\cP)$  is a proposition involving the orbits of a dynamics which may be true, wrong or undecidable. This defines a function:
\[\cP: C^r(M,M)\to \{\mathrm{true, wrong, undecidable}\}\; .\]
We would like to study properties which, roughly speaking, occur for all systems but some in a negligible set in some rigorous sense. Such properties will be said to be \emph{typical}. A similar way of description already occurred in singularity theory of differentiable mappings: Whitney \cite{Wh34} showed that every closed set of a manifold $M$ is the zero set of a function. That is why Thom and then Mather restricted their study of singular set to  differentiable maps in an open and dense set. 

In this article we will see that  \emph{typically} the behavior of a dynamical system is much more complex than expected, for properties involving the growth of the number of periodic points, the value of entropy, or a recent notion of emergence. This complexity will regard mostly localizable properties.
\begin{definition}
A property $(\cP)$ on $\Diff^r(M, M)$ is \emph{localizable}, if the following condition holds true:
for every open set $U\subset M$, $N\ge 0$  and $f\in  \Diff^r(M)$ such that $f^N|_U=id|_U$, there exists a $C^r$-perturbation $\tilde f$ of $f$ which satisfies $(\cP)$. 

A property $(\cP)$ on $\Diff^r(M)$ is \emph{openly localizable}, if the perturbation $\tilde f$ can be chosen in a $C^r$-open set.

A property $(\cP)$ is \emph{genericaly localizable}, if it is implied by the conjunction of a countable number of openly localizable properties $\{(\cP_i): {i\ge 0}\}$: $$\bigwedge_i (\cP_i)\Rightarrow (\cP)\; .$$ 
 \end{definition}
This definition is inspired by the following result:
\begin{theorem}[Turaev \cite{Tu15}]\label{Turaev}
For every surface $M$ and $r\ge 2$, there is a non-empty open set  $\cal N^r\subset \Diff^r(M)$ such that for a dense set $\cal D\subset \cal N^r$, for every $f\in \cal D$, there exist $U\subset M$ and $N\ge 1$ such that $f^N|_U=id|_U$.
\end{theorem}
An immediate consequence of this result is:
\begin{corollary}\label{Turaev consequence} Among surface $C^r$-diffeomorphisms, $r\ge 2$, a generically localizable property $(\cP)$  is locally generic: there exists $\cal R\subset \cal N^r$ topologically generic (as defined in the sequel) such that every $f\in \cal R$ satisfies $(\cP)$.
\end{corollary}
 \medskip 
 
We will discuss in the next sections about the following localizable properties: \begin{enumerate}
\item displaying a fast growth of the number of periodic points (which is generically localizable),
\item having positive metric entropy (which is localizable),
\item having a  high emergence (which is generically localizable).
\end{enumerate}
As a matter of fact, these properties are valid on sets which are typical for different notions.

\section{Typicality}
In this section we recall and develop some notions of typicality.  
This will enable better interpretations of Problem \ref{problem2} and also Conjectures \ref{positive entropy conjecture} and \ref{positive order emergence conjecture} which are part of the dictionary given \cref{dico entropy and emergence }.

For the sake of simplicity, let us assume that $M$ is compact. Let $1\le r\le \infty$ and let $\cal U$ be a non-empty,  open subset of  $C^r(M,M)$ endowed with its canonical topology. 

\subsection{Topological notions of typicallity}
 From the topological viewpoint, there are the following notions of typicality:
\begin{definition} 
A property  is \emph{open and dense} in $\cal U$ if it  holds true at an open and dense subset of $\cal U$. 
A property  is \emph{topologically generic} in $\cal U$  if it  holds true at a countable intersection of open and dense subsets of $\cal U$
.  A property is \emph{dense} in $\cal U$ if it holds true at a dense subset of $\cal U$\end{definition}

As $\cal U$  is a Baire space (a countable intersection of open and dense sets is a dense set), the above notions of typicality are linked  as follows, where $(\neg\cP)$ is the negation of $(\cP)$:
\[\left\{\begin{array}{c}
(\cP) \text{ is open and dense} \Rightarrow (\cP) \text{ is topologically generic}  \Rightarrow (\cP) \text{ is dense}\\
(\neg\cP)  \text{ is dense}\Rightarrow (\cP)\text{ is \emph{not} open and dense}\\
 (\neg\cP) \text{ is topologically generic} \Rightarrow (\cP) \text{ is \emph{not} topologically generic}\end{array}\right.
\]

The main failure of density as a notion of typicality  is that the intersection of two dense sets might be empty (and so not dense), whereas it is reasonable to ask for the typicality of the conjunction of two typical properties.

On the other hand, a finite intersection of open and dense sets is still open and dense, and a countable intersection of topologically generic sets is topologically generic as well.  However, these notions of typicality display others problems.  The notion of open and dense property as typicality is very natural, but in practical it is often too restrictive: for instance having all periodic points with eigenvalues $\neq 1$ is not an open and dense property (this has been proved first in \cite{Newhouse}, and is also a consequence of  \cref{Turaev}). 

The main problem with the notion of topological genericity is that such subsets might be negligible in the measure theoretical viewpoint. Indeed, there are topologically generic subsets of $\R$ whose Lebesgue measure is zero. That is why Kolmogorov \cite{Ko57}, motivated by von Neumann \cite{vN43}, proposed a notion of typicality which involves probability. 

\subsection{Probabilitic notions of typicality}
In his ICM plenary talk \cite{Ko57}, Kolmogorov briefly proposed to consider deformation of a system $f_0$ along a parameter family $(f_a)_a$ such that 
\[f_a(x)= f_0(x)+a \phi(x, a)\]
and wonder if a property persists for Lebesque almost every parameter $a$ small.  Kolmogorov was proposing an initial working frame with $\phi$ analytic, but with perspective in finite regularity $C^r$, $r<\infty$.  The differentiable counterpart of this notion should involve $C^r$-families of mappings $(f_a)_a$ passing through the map $f_0$. However, the situation is here much less rigid, and one must avoid the family to lie in an exceptional region. That is why one need to re-interpret this notion. There were several interpretations of this notion of typicality \cite{IL99, KH07}. Let us chose the following, involving the unit ball $B^d$ of $\R^d$:
\begin{definition}\label{Ktypical}
A property $(\cP)$ is \emph{$d$-Kolmogorov typical} in  $\cal U$, if there exists a Baire residual set $\cal R \subset \{(f_a)_a\in C^r(B^d\times M;M): f_a\in \cal U,\; \forall a\in B^d\}$ such that for every $(f_a)_{a\in B^d}\in \cal R$, for Lebesgue a.e. parameter  $a\in B^d$, the dynamics $f_a$ satisfies property $(\cP)$.
\end{definition}
Note that  $0$-Kolmogorov typicality is the topological genericity. For $d>0$, $d$-Kolmogorov typicality implies the density of a phenomena, but not necessarily its topological genericity.

We notice that a $d$-$C^r$-Kolmogorov typical subset contains Lebesgue a.e. points of (unaccountably many) submanifolds of $C^r(M,M)$ of dimension $d$.  On the other hand, one can construct (artifical) $d$-$C^r$-Kolmogorov typical sets which intersects any $d+1$-submanifold of $C^r(M,M)$ at a set of Lebesgue measure 0 \cite{IL99}. The latter point is the main defect of this notion. To overcome this, a way  is to introduce:
\begin{definition}\label{Ktypicalgen}
A property $(\cP)$ is \emph{eventually Kolmogorov typical} in  $\cal U$, if it is $d$-Kolmogorov typical for every $d\ge 0$ sufficiently large.
\end{definition}

Another similar notion, which is not weaker nor stronger is the following:
\begin{definition}\label{Defabundantl}
A property $(\cP)$ is \emph{$d$-$C^r$-abundant} in $C^r(M;M)$, if there exist $d\ge 1$ and a non-empty open subset $\cal V \subset  C^r(B^d\times M;M)$ such that for every $(f_a)_{a\in B^d}\in \cal V$, there is a parameter set $E\subset B^d$ of positive Lebesgue measure such that for every $a\in E$, the map $f_a$ satisfies property $(\cP)$.
\end{definition}
Note that there exist  properties which are abundant and whose contrary is also abundant. This notion is related to the Kolmogorov typicality as follows:
\[ \neg(\cP) \text{ is $d$-$C^r$-Kolmogorov typical} \Rightarrow (\cP) \text{ is \emph{not} $d$-$C^r$-abundant}\]
As the relevance of abundance is incontestable, so is the notion of  Komogorov typicality.

In the next years I will be working on the following conjecture:
\begin{conjecture}\label{Conjecture}
If the dimension of $M$ is sufficiently large, for every $2\le r\le \infty$, every property which is generically localizable in $\Diff^r(M)$ is $d$-Kolmogorov typical in some open set of $\Diff^r(M)$, for every $d\ge 0$.
\end{conjecture}

For the sake of completeness, we shall recall the notion of prevalence \cite{HSY92}, which is another probabilistic notion of typicality. 
A  subset of a vector space $V$  is \emph{prevalent} if its complement is shy. A set $S$ is \emph{shy} if there exists a measure $\mu$ on $V$ which has positive and finite value on a compact set and such that $\mu(S+v)=0$ for every $v\in V$. An important case is when $\mu$ coincides with the  Lebesgue measure on a finite dimensional set $E\subset B$. Then a set is shy if at every $p\in S$, the affine line $p+E$ intersects $S$ at a subset of Lebesgue measure $0$.

A substantial problem of the notion of shy \& prevalence is that the parameter space $C^r(M,M)$ does not display a canonical vector structure in general. Of course  it is a Frechetic manifold and so it is locally modeled by a Frechet space, however a coordinate change of $C^r(M,M)$ nor of $M$ does not preserve this structure. In particular, prevalence has no reason to be invariant by these coordinate changes.

\section{Growth of the number of periodic points}
A natural dynamically defined set is the one of its $n$-periodic points $$Per_n f:= \{x\in M: f^n(x)=x\}\; .$$
 To study its cardinality, we consider also the subset $Per_n^0 f\subset Per_n f$ of its isolated points. We notice that the cardinality of  $Per^0_n f$ is invariant by conjugacy. Hence it is natural to study the growth of this cardinality with $n$.

Clearly, if $f$ is a polynomial map, the cardinality $Per^0_n f$ is bounded by the degree of $f^n$, which grows at most exponentially fast \cite{Milnor-Friedman,DNT16}.

The first study in the $C^\infty$-case goes back to Artin and Mazur \cite{AM65} who proved the existence of a dense set $\mathcal D$ in $\Diff^r(M)$, $r\le \infty$, so that for every $f\in \mathcal D$, the number $Card\, Per^0_n f$  grows  at most  exponentially, i.e. :
\begin{equation}\tag{A-M} \limsup_{n\to \infty}  \frac1n \log Card\, Per^0_n f<\infty 
\; .\end{equation}

This leads Smale \cite{Sm} and Bowen \cite{Bo78} to wonder whether (A-M) diffeomorphisms are generic. Finally Arnold asked the following problems:
\begin{problem}[Smale 1967, Bowen  1978, Arnold  Pb. 1989-2 \cite{Ar00}] \label{problem1} Can the number of fixed points of the $n^{th}$ iteration of a topologically generic infinitely smooth self-mapping of a compact manifold grow, as $n$ increases, faster than any prescribed sequence $(a_n)_n$ (for some subsequences of time values n)? 
\end{problem}
This problem enjoys a long tradition.  In dimension 1, Martens-de Melo-van Strien \cite{MdMvS92} showed that for every $\infty\ge r\ge 2$, for an open and dense set\footnote{whose complement is the infinite codimentional manifold formed by maps with at least one flat critical point.} of $C^r$-maps, the number of periodic points grows at most exponentially fast.

In higher dimension, Kaloshin \cite{K99} proved the existence of a dense set  $\mathcal D$ in $\Diff^r(M)$, $r< \infty$, such that the cardinality of $Per_n f\supset Per_n^0f$ is finite and grows  at most  exponentially fast. On the other hand, he proved in \cite{K00} that whenever $2\le r<\infty$ and $\dim M\ge 2$,  \emph{a locally topologically generic diffeomorphism displays a fast growth of the number of periodic points}: 
there exists an open set $\cal U\subset\Diff^r(M)$,  such that for any sequence of integers $(a_n)_n$, a topologically generic $f\in \cal U$ satisfies:
\begin{equation}\tag{$\star$}
\limsup_{n\to \infty} \frac{Card\, P_n^0(f)}{a_n}=\infty\; .\end{equation} 
\begin{proposition}\label{growth per point}
To display a fast growth on the number of periodic points is a generically localizable property in $\Diff^r(M)$, for every $r\ge 1$.
\end{proposition}
\begin{proof}
Let $(a_j)_j\in \N^\N$. Let $(\cP_j)$ be the property that $f$ has $\ge a_j$-saddle points of period $\ge j$.   
 
We observe that $(\cP_j)$ is an openly localizable property. Indeed, 
let  $ U \subset M$  be an open subset and $N$ an integer such that $f^N|_U=id|_U$. 
We may take $N$ and $U$ smaller such that $f^k(U)$ does not intersect $U$, and then do a small perturbation $\tilde f$ supported by $U$ which satisfies $(\cP_j)$. The property is open by hyperbolic continuation of the saddle points. Obviously the conjunction $\bigwedge_j (\cP_j)$ implies the fast growth of the number of periodic points. 
\end{proof}

\cref{Turaev consequence} of  Turaev's Theorem  and \cref{growth per point} implies that  Kaloshin's theorem is also valid in the space of  $C^\infty$-surface diffeomorphisms. In dimension $\ge 3$, this result has been extended to the $C^1$-case and to the $C^\infty$-case in \cite{BerPer2018,AST2018}. 
 In summary, we have the following description:
\begin{theorem}\label{theoA}
Let $\infty\ge r\ge 2$ and  let $M$ be a compact manifold of dimension $n$.
 \begin{itemize}
\item
If $n=1$, Property $(AM)$ is satisfied by an open and dense set of $C^r$-self-mappings.
\item
If $n\ge 2$, there exists a (non-empty) open set $\cal U\subset\Diff^r(M)$ so that given any sequence $(a_n)_n$  of integers, a topologically generic $f$ in $\cal U$ satisfies $(\star)$.
\end{itemize}
 \end{theorem}
 The conservative counterpart of this result was established in \cite{KS06}. For a stronger sense of typicality, it has been asked \footnote{Many other Arnold's problems are related to this question \cite[1994-47, 1994-48, 1992-14]{Ar00}.}:
\begin{problem}[Arnold 1992-13 \cite{Ar00}]\label{problem2} Prove that a typical, smooth, self-map $f$ of a compact manifold satisfies that $(Card\, Per_n f)_n$ grows  at most exponentially fast.  
\end{problem}

To provide a positive answer to this Problem, Hunt and Kaloshin \cite{KH07,KH072} used a method described in \cite{GHK06} to show that for $\infty\ge r>1$, a \emph{prevalent}  $C^r$-diffeomorphism (for a chosen vector structure) satisfies:
\begin{equation}\tag{$\Diamond$}
\limsup_{\infty} \frac{\log P_n(f)}{n^{1+\delta}}= 0,\quad \forall \delta>0\; .\end{equation}

However the latter did not completely solve Arnold's Problem \ref{problem2} in particular because the notion of prevalence is \emph{a priori}  independent to the notion of typicality initially meant by Arnold. Indeed this problem was formulated for typicality\footnote{See explanation below problem 1.1.5 in \cite{ KH07}.} in the sense of definition \ref{Ktypical}. Actually following this notion, the  answer to  Arnold's Problem \ref{problem2} is negative in the finitely smooth case:

\begin{theo}[\cite{BerPer2018}]
\label{theoB}
Let $1\le r< \infty$, $0\le d<\infty$ and let $M$ be a manifold of dimension $\ge 2$.  

Then there exists a (non-empty) open set $\hat {\cal U}$ of $C^r$-families $(f_a)_{a\in B^d}$ of $C^r$-self-mappings $f_a$ of $M$ such that, for any sequence of integers $(a_n)_n$, a topologically generic $(f_a)_a\in \hat {\cal U}$ consists of maps 
$f_a$ satisfying $(\star)$ for every $a\in B^d$.
Moreover if $\dim \, M\ge 3$, the maps $f_a$ are diffeomorphisms.
\end{theo}
The symplectic counterpart of the above result has been proved previously in \cite{As16}. Both confirm Conjecture \ref{Conjecture}.

\section{Entropy}

In his revolutionary work \cite{poincare1892methodes}, Poincaré showed the existence of an open class of systems which display  an infinite, invariant,  set $E$ 
in which all the points have pairwise different asymptotic behaviors. His example was given by the so-called restricted 3-body problem in celestial mechanic. Mathematically, it is a Hamiltonian, analytic maps of a surface. 
Nowadays, the dynamics of this example (as most of the symplectic dynamics) is poorly understood. 
For instance, we do not know if $E$ has zero Lebesgue measure for the Poincaré example.

The concept of ergodicity -- introduced by Boltzmann -- may enable to understand most orbits of such systems. The idea is that the statistical behavior of the orbit may be the same for ``most'' of the point in $E$. The statistical behavior of the orbit of a point $x$ is given by the sequence  of measures $\mathsf e_n(x)$ equidistributed on the  first iterates up to time $n$. These are  called the \emph{empirical measures}, and are defined as follow (when the time is discrete):
\[\mathsf e_n(x) =\frac1n \sum_{0\le k< n} \delta_{f^k(x)}\; .\]
A probability measure $\mu$ is \emph{ergodic} if for $\mu$-a.e. every point $x$, the sequence $(\mathsf e_n(x))_n$ converges to $\mu$. In other words, the time average of the orbit of $\mu$-a.e. point is equal to the space average.
The first example of differentiable realization of such systems was done first for billiards \cite{Ha98} and then for geodesic flow of negatively curved compact manifolds by  Hedlund  and  Hopf  \cite{He39,Ho39}. 

To quantify the sensitivity to the initial conditions, we consider the sequence of distances:\label{def dn}
\[d_n(x,y):= \max_{0\le i<n} d(f^i(x), f^i(y))\; .\]
The concept of entropy was introduced by Kolmogorov in dynamics, and reformulated as follows by Katok \cite{Ka80}. 
\begin{definition}[Entropy]\label{def metric entropy} 
Given an ergodic measure $\mu$, let $\cal N^\mu_n(\epsilon)$ be the minimal number of $d_n$-balls of radius $\epsilon>0$ whose union has $\mu$-measure $\ge 1-\epsilon$. The entropy $h_\mu$ of $\mu$  is:
\[h_\mu:= \lim_{\epsilon \to 0}\lim_{n\to \infty} \frac1n \log \cal N^\mu_n(\epsilon)\; .\]
\end{definition} 
Hence entropy measures the complexity of the asymptotic behavior.  It turns out that ergodic measures of positive entropy are well modeled.
By Austin's theorem \cite{A18}, every ergodic automorphism of positive entropy is isomorphic to the direct product of a Bernoulli shift and a 
transformation of arbitrarily small entropy. 
 
After the work of Anosov and Sinai, it has been shown that the Liouville measure under the action of a  geodesic flow of a negatively curved surface has positive entropy. This example has been generalized to discrete time, by Anosov and then Smale \cite{Sm} to reach the notion of uniformly hyperbolic dynamics. Such systems are such that their non-wandering set is locally maximal and hyperbolic. An invariant set $K$ is \emph{hyperbolic}  if there is a continuous splitting $E^s\oplus E^u=TM|_K$ by two invariant bundles $E^s$ and $E^u$ which are respectivelly contracted and expanded by the differential of the dynamics. By the work of Sinai, Ruelle and Bowen \cite{BR75}, we know that every uniformly hyperbolic smooth dynamics 
displays a finite number of ergodic probability measures whose basins cover all the manifold modulo a set of Lebesgue zero.  We recall that the \emph{basin} $B_\mu$ of an ergodic probability measures $\mu$ is
\[B_\mu:= \{x\in M: \mathsf e_n(x)\rightharpoonup \mu\}\; .\]
 In this extend, a special interest is given ergodic measure which are \emph{physical}: those whose basin has positive Lebesgue measure.

\subsection{The positive entropy conjecture}
%
A diffeomorphism is said to have \emph{positive (metric)  entropy} if it leaves invariant an ergodic, physical probability with positive entropy.   Here is one of the most important and difficult question in this field:
\begin{conjecture} \label{positive entropy conjecture}
A typical symplectic, surface diffeomorphisms has positive entropy.
\end{conjecture}
Here ``typical'' may be view as occurring in some natural examples (from physical science): for instance the Poincaré example, or as conjectured by Sinai \cite[P. 144]{Si94}, the Chirikov standard maps. This may be seen also following one of the notions of typicality of the first section. 

Herman proposed the following (weak) version of the positive entropy conjecture:
\begin{conjecture}[\cite{He98}]
There are conservative, infinitely smooth perturbations of the identity of the $2$-disk with positive metric entropy.
\end{conjecture}
Note that this conjecture implies that having positive entropy is a localizable property. 

In the same paper, he also wondered if the set of surface, conservative maps of the disk with  positive metric entropy are dense (in the $C^\infty$-topology), or if its interior is non-empty. Herman's conjecture has been solved recently:
\begin{theorem}[\cite{BT19}] 
In the open set of conservative $C^\infty$-surface maps which display an elliptic periodic point, those with positive metric entropy form a dense subset.
\end{theorem}
A map which  robustly does not display any periodic elliptic point is called \emph{weakly stable}. The conservative counterpart of the Ma\~ ne conjecture \cite{Ma82} states that weakly stable maps are uniformly hyperbolic. As a conservative, uniformly hyperbolic map has positive entropy, this theorem implies, \emph{up to the latter  conjecture}, the density of the positive metric entropy. This theorem was proved by showing that a dynamics constructed by surgery from a uniformly hyperbolic system (similar to Przytycki \cite{Pr82} example), appears after perturbation and renormalization. We know that these examples do not appear among entire maps \cite{Us80}. Also  these examples are not locally generic nor typical in the sense of Kolmogorov \cite{He83}.

Hence to solve Conjecture \ref{positive entropy conjecture} (if it is correct!), one needs to find new models.  In order to do so, Yoccoz' strong regularity program proposes \cite{Yolecture1,BY19} to look for geometrical and combinatorial definition of invariant sets which display an ergodic, physical measure with positive entropy. The idea of this program  is to work on paradigmatic and abundant examples of (not necessarily conservative) smooth surface maps, such that the dimension of the measure is higher and higher. The \emph{dimension of a measure} $\mu$ is the infimum of the Hausdorff dimension of a set $E$ such that $\mu(E)=1$.

The first paradigmatic examples were Yoccoz' strongly regular quadratic maps \cite{Y,Y19} which allowed him to give an alternative proof of Jakobson theorem \cite{Ja81}. In \cite{PY09}, the strongly regular horseshoes were introduced and shown to be abundant. In \cite{berhen}, strongly regular Hénon-like attractors were introduced and shown to be abundant. This gives an alternative proof of Benedicks-Carleson theorem \cite{BC2} and solves the second step  of Yoccoz's program as stated during his first lecture at Coll\`ege de France (and the main source of inspiration of his last, unfinished lecture \cite{Yo16}).  In \cite{BY19}, a few conjectures are stated to develop this theory.

\subsection{Description of measure with positive entropy}
Pesin theory describes the ergodic measures of positive entropy of a (possibly not conservative) diffeomorphism of class $C^r$, with $r>1$. 

We recall that by Oseledets' theorem, for any ergodic measure $\mu$,  there exists a $Df$-invariant splitting  $T_xM=
E^s(x)\oplus E^c(x)\oplus E^u(x)$ depending measurably of $\mu$ at a.e. $x$ and such that:
\[ \lim_{n\to \pm \infty}  \frac 1n \log \|D_xf^n| E^s\|<0\quad 
\lim_{n\to \pm \infty} \frac 1n \log \|D_xf^n|E^c\|=0 \qand 
\lim_{n\to \pm\infty} \frac 1n \log \|D_xf^n|E^u\| >0\; .\]
Furthermore, with $u= \dim\, E^u$, there exists a flag $0\subset E^1\subset E^2\subset  \cdots\subset  E^{u-1}\subset E^u$ 
associated to values $0<\lambda_1\le \cdots \le \lambda_{u}$ such that any vector of $E_{i+1} \setminus E_i $ is asymptotically expanded by a factor $\exp(\lambda_i)$.
The numbers $0<\lambda_1\le \cdots \le \lambda_{u}$ are called the \emph{positive Lyapunov exponents} and their sum is denoted by $\Sigma_+$. 
\begin{theorem}[Ruelle inequality \cite{Ru78}]
\label{Ruelle ineq} It holds: 
$h_\nu \le \Sigma_+ \le \dim M\cdot \int \log\| Df \| d\nu$.
\end{theorem}
Hence a measure with positive entropy has necessarily a positive Lyapunov exponent.  
By the Pesin theorem \cite{Pe76, FHY83}, for $\mu$-a.e. $x$,     there exist immersed manifolds $W^u(x)$ and $W^s(x)$ tangent to respectively $E^u(x)$ and $E^s(x)$ such that 
\[W^s(x; f):=\left\{y\in M: \overline {\lim_{n\to +\infty}} \frac 1n \log d(f^n(x),f^n(y))<0\right\} \qand 
 W^u(x; f):=  W^s(x; f^{-1})\; .\]
Similarly for every $i$,  the following is an immersed manifold $W^i(x)\subset W^u$ tangent to $E^i(x)$:
\[W^i(x; f):=\left\{y\in M: \overline {\lim_{n\to \infty}} \frac 1n \log d(f^{-n}(x),f^{-i}(y))\le -\lambda_i\right\}\; .\]
Given $\epsilon>0$, let $W^i_\epsilon(x; f)$ be the connected component of $x$ in the intersection of $W^i_\epsilon(x; f)$ with the $\epsilon$-neighborhood of $x$.  These subsets vary measurably in function of the point $x$, and so for every $i$,  there is a compact set $T_i$ such that $(W^i_\epsilon (x))_{x\in T_i}$ is a continuous disjoint family of $C^2$-subsets, whose union has positive $\mu$-measure. Then there exists a family of measure $\mu_x^i$ and a measure $\nu_i$ on $T_i$ such that for every $A\subset \bigcup_{x\in T_i} W^i_\epsilon (x)$, it holds:
\[\mu(A)= \int_{T_i} \mu_x^i(A\cap W^i_\epsilon (x)) d\nu_i(x)\; .\]

The measure $\mu$ is SRB (Sinai-Ruelle-Bowen) if $\mu_x^{u}$ is absolutely continuous w.r.t. the Lebesgue measure of $W^u_\epsilon(x; f)$ for $\mu$-a.e. $x$. If the measure is furthermore \emph{hyperbolic} (i.e. $E^c=0$), then  the measure $\mu$ is physical. Also by the Ledrappier-Strelcyn Theorem,  its entropy is positive and equal to $h_\mu= \sum_i (dim(E_i)-\dim(E_{i-1}))\lambda_i$. 

More generaly, Young defined the dimension of the measure $\mu_i^x$ as:
\[d_i:= \lim_{\epsilon \to 0}\frac{\log\mu_i^x(W^i_\epsilon(x; f(x))}{\log \epsilon}\; .\]
This is invariant by the dynamics, and so by ergodicity, this does not depend on $x$. We have:
\begin{theorem}[Ledrappier-Young entropy formula \cite{LYII85}]\label{Ledrappier-Young entropy formula} 
$h_\mu = \int \sum_i \lambda_i (d_i-d_{-i-1})d\mu(x). $
\end{theorem}
 The proof of this theorem sheds light on a certain fractal structure of the support of the measure restricted to the submanifold $W^i_\epsilon(x;f)$. This will be a source of inspiration for the next section.

\subsection{Measure of maximal entropy}
Given a $C^1$-map $f$ of a compact manifold. For every $n\ge 0$, let  $d_n$ be the  distance   defined \cpageref{def dn}. Let $\cal N_n(\epsilon)$ be the minimal number of $d_n$-balls of radius $\epsilon>0$ covering $M$. 
\begin{definition}\label{def topological entropy} The topological entropy of $f$ is:
\[h_{\mathrm{top}}:= \lim_{\epsilon \to 0}\lim_{n\to \infty} \frac1n \log \cal N_n(\epsilon)\; .\]\end{definition}
Let us recall the following:
\begin{theorem}[Variational principle \cite{Di71}]\label{varia entropy}
The topological entropy is the supremum of the metric entropy among all the ergodic probability measures.
\end{theorem}
For many transitive sets \cite{Ma83, Lyu83, BLS, Pa64, berent, BCS18}, the maximal entropy measure exists and is unique. Interestingly,  the  maximal entropy measure $\mu_h$ is often equi-distributed on the periodic points:
\[\frac1{\mathrm{ Card\, Fix}(f^n)} \sum_{x\in Fix(f^n)} \delta_x\to \mu_h\; .\]
This property has been shown for any rational function (of positive degree) of the Riemanian sphere by \cite{Ma83, Lyu83}, any complex polynomial automorphisms (of positive dynamical degree) by \cite{BLS}, for strongly regular Hénon map in \cite{berent}, and recently \cite{Bu17} for smooth  diffeomorphisms whose  ergodic measures are all hyperbolic with Lyapunov exponents uniformly  non-zero.

\section{Emergence}
Metric entropy describes very well the complexity of a dynamics of a.e. point of  an  ergodic probability measure with positive entropy. This would be sufficient if Bolzman ergodic hypothesis would be correct.
In modern language \cite{BK32}, this hypothesis states that typical Hamiltonian systems are ergodic on [$\Leb$. almost] every [component] of energy level.  However, this hypothesis turned out to be wrong after a development of the KAM theory \cite{Mo62}. This theory states that any perturbation of some integrable systems displays an invariant set of positive Lebesgue measure which is laminated by invariant Lagragian tori on which the dynamics acts as an irrational rotation.  

To get rid of this obstruction, Smale proposed to simplify this study, by considering differentiable systems of surface which do not leave invariant a volume form. At some point he conjectured that an open and dense set of such systems should be uniformly hyperbolic. As we saw above, by \cite{BR75}, this would have implied that an open and dense set of differentiable dynamical systems satisfy the following property: 
there exists a finite number of ergodic probability measures  such that the union of their basins contains  Lebesgue a.e. point of the  manifold. 
This scenario (and Smale's conjecture) was disproved by his student Newhouse \cite{Newhouse} who discovered that the following phenomenon appears topologically generically in an open set of $C^r$-surface diffeomorphisms, $2\le r\le \infty$:

\begin{defi}[Newhouse phenomenon]
There exists an infinite, transitive hyperbolic compact set which is accumulated by infinitely many sinks  or totally elliptic periodic points.
\end{defi}

Interestingly this phenomenon appears in many fields of dynamical systems:
\begin{center}

\medskip
\begin{tabular}{|r|l|}
  \hline
  & Newhouse phenomenon is topologically generic in an  open set of:
\\  \hline
Newhouse \cite{
Newhouse}&$C^r$-surface diffeomorphisms for every $r\ge 2$.\\
\hline
  Duarte \cite{Duarte99
  } & $C^r$-surface, conservative diffeomorphisms for every $r\ge 2$.\\
  \hline
 Buzzard \cite{Bu97} & Polynomial automorphisms of $\C^2$ of some large degree. \\
   \hline
Bonatti-Diaz \cite{BD99} & $C^1$-diffeomorphisms of manifold of dimension $\ge 3$.\\
  \hline
 Biebler \cite{biebler2016persistent} & Polynomial automorphisms of $\C^3$ of any degree $\ge 2$.\\
  \hline
\end{tabular}
\end{center}
\medskip

We will see below that the set of invariant probability measures is huge. In this extend, this phenomenon is very much in opposition with 
Bolzman ergodic hypothesis. 

Several conjectures from the 90's \cite{TLY, PT93, PS95, Pa00, Pa05, Pa08} gave some (philosophical) hope that Newhouse phenomenon should be negligible in some sense of Kolmogorov. But this turns out to be \emph{not} the case in finite differentiability  $1\le r<\infty$ and dimension $\ge 2$: 
\begin{theorem}[\cite{BE15, Be16}] 
Let $\infty> r\ge 1$ and $0\le d<\infty$, let $M$ be a manifold of dimension $\ge 2$.  Then there exists a (non-empty) open set $\hat {\cal U}$ of $C^r$-families $(f_a)_{a\in B^d}$ of $C^r$-self-mappings $f_a$ of $M$ so that, a topologically generic $(f_a)_a\in \hat {\cal U}$ satisfies that for every $\|a\|\le 1$,  the dynamics $f_a$ displays Newhouse phenomenon.

Moreover if $dim \, M\ge 3$, we can chose $\hat {\cal U}$ to be formed by families of diffeomorphisms.
 \end{theorem}

Still one can argue that in Newhouse phenomenon, the Lebesgue measure of the basins of the sinks of period $n$ might decrease very fast, and a typical mapping (including those displaying Newhouse phenomenon) might have their statistical behavior which is well approximated by a finite number of probability measures. Let us study this possibility.
\subsection{Quantitative approach} 

The aim of the notion of emergence is to quantify how far a system is to be ergodic. Let us start with the abstract setting where the dynamics $f$ is a continuous map of a compact metric space $X$ with well defined box dimension and endowed with a reference probability  measure $\mu$. The important cases will be when $X$ is a compact manifold and $\mu$ a Lebesgue measure.  We endow the space of probability measures $\mathcal M(X)$ of $X$ with the following distance:

\begin{definition}[Kantorovich-Rubinstein distance ${\mathbf d}$] For every $\nu_1, \nu_2\in \cal M(X)$, put:
\begin{equation*}{\bf d}(\nu_1, \nu_2)=\sup_{\phi\in Lip^1(X)} \int \phi \, d(\nu_1-\nu_2),  \end{equation*}
with $Lip^1(X)$ the space of real functions on $X$ which are $1$-Lipschitz.
\end{definition}
We recall that the statistical behavior of the system is given by the empirical measures $\mathsf e_n(x) =\frac1n \sum_{0\le k< n} \delta_{f^k(x)}$, for $x\in X$  and $n\ge 0$. Roughly speaking, the {emergence} ${ \cE_\leb(\epsilon)}$ at scale $\epsilon$ is the number of probability measures needed to approximate the statistical behavior of the system with precision~$\epsilon$:
\begin{definition}[Metric emergence \cite{Be16}]
The \emph{ metric emergence} $\cE_\mu (f)$ of $f$ at scale $\epsilon>0$ is the minimal number $N$ of probability measures $\{\nu_i\}_{i\le N}$ so that 
\begin{equation}\tag{$\diamondsuit$}
\limsup_n \int_M   {\bf d}(\mathsf e_n(x), \{\nu_i\}_{i\le N})d\mu <\epsilon\; .\end{equation}
\end{definition}
We recall that if $f$ leaves invariant $\mu$, then by the Birkhoff ergodic theorem, for $\mu$-a.e. $x$, the sequence $(\mathsf e_n(x))_n$ converges to a measure $\mathsf e(x)$. Then $(\diamondsuit)$ is equivalent to:
\[\int_M {\bf d}(\mathsf e(x), \{\nu_i\}_{i\le N})d\mu<\epsilon  \; .\]
Let us recall the following examples:
\begin{exem}[\cite{Be16,BB19}] Let $X$ be a compact manifold of dimension $d$ and $\mu=\leb$:

\begin{itemize} 
\item If $\leb$ is ergodic for $f$ (e.g. $f$ is an Anosov map), then $\cE_\leb(f)=1$. 
\item A dynamics displaying Newhouse phenomenon has not finite emergence. 
\item The identity satisfies  $\lim_{\epsilon\to 0} \frac{\log \cE_\leb (id)}{-\log \epsilon}=d$.
\item If $M$ is symplectic and $f$ displays KAM phenomenon, then $ \frac{\log \cE_\leb (f)}{-\log \epsilon}\ge d/2$.
\end{itemize} 
\end{exem}
 All the latter bounds from below are at most polynomial ($\epsilon^{-d}$). We proposed the following notion to describe dynamics whose statistical behavior is  poorly described by finitely many probability measures:
  \begin{definition}
A  differential dynamics is with \emph{high emergence} if:
\[\tag{Super Polynomial Emergence}{ \limsup_{\epsilon \to 0} \frac{\log \cE_\leb(\epsilon)}{-\log \epsilon}=\infty}\; .\]
 \end{definition}

 The first aim of our program on Emergence is to prove the following conjecture:
\begin{conjecture}[\cite{Be16}]\label{positive order emergence conjecture}  Dynamics of high emergence are  typical in many senses and contexts.\end{conjecture} 
The second aim of our program on emergence is to build a theory describing differentiable dynamics with high emergence.
\subsection{Results on typicality of dynamics with high emergence}
There are basically two ways to obtain high emergence:
\begin{enumerate}
\item[(\emph a)] For every $x$ in a set of positive Lebesgue measure, the sequence of empirical measures $(\mathsf e_n(x))_n$ does not converge and its cluster values form an infinite dimensional set,
\item[(\emph b)] At a.e. $x$  the empirical measures $(\mathsf e_n(x))_n$ converges to a measure $\mathsf e(x)$, but the ergodic decomposition $\mathsf e_*\leb$ is poorly approximated by a finitely supported measure on $\cal M(M)$. \end{enumerate}

In our bounds on emergence, the scale will be often the following:
\begin{definition}[Order] Given a function $\phi: (0, \infty)\to (0, \infty)$, its \emph{order} is:
$${ \mathcal {O}\phi}:=\limsup_{\epsilon\to  0} \frac{\log \log \mathcal \phi(\epsilon)}{-\log \epsilon}.$$
\end{definition}
Let us observe that the emergence is high if the order of the emergence  is positive:
 \[ \mathcal {O}\cE_\leb:=\limsup_{\epsilon\to  0} \frac{\log \log \mathcal \cE_\leb(\epsilon)}{-\log \epsilon}>0\Rightarrow 
  \lim_{\epsilon\to 0} \frac{\log \mathcal \cE_\leb(\epsilon)}{-\log \epsilon}=\infty\; .\]

We observe that the metric emergence is at most the covering number\footnote{The \emph{covering number} at scale $\epsilon$ of a metric space is the least number of $\epsilon$-balls which cover it.} of the space of probability measures, which is bounded by the following:
 \begin{theorem}[Bolley-Guillin-Villani \cite{BGV}, Kloeckner \cite{Kl12},  Berger-Bochi-Peyre \cite{BB19}]\label{analogue Ruelle} If $X$ is a compact space with well defined box dimension, then the covering number $H(\epsilon)$ of ($\mathcal M(X), {\bf d})$ has order $\mathcal O H=\dim X\; .$
\end{theorem}
Hence we obtain the following counterpart  of the Ruelle inequality formula for emergence:
\begin{corollary}
If $f$ is a continuous map of a compact metric space $X$ with  box dimension $d$, and $\mu$ a probability measure on $X$, it holds:
\[\cal O \cE_\mu(f)\le d.\]
\end{corollary}

This maximum is attained in generic differentiable examples:
\begin{theorem}[\cite{BB19}]\label{typicality conserva emer}
A topologically generic $f\in \Diff^\infty_\Leb(M^2)$ which displays an elliptic periodic point has its metric emergence of maximal order $\dim M^2=2$:
$$\mathcal {OE}_\leb= \dim M^2\; .$$
\end{theorem}
Let us recall that maps with robustly no elliptic points (which are called \emph{weakly stable}) are conjecturally uniformly hyperbolic (Ma\~ né-Like  conjecture \cite{Ma82}) and ergodic. Hence conjecturally, the emergence of a generic conservative maps is either minimal\footnote{The minimal emergence equal 1, in this conservative setting, this means that the Lebesgue measure is ergodic.} or maximal.

The latter theorem, and the next one were proved using that maximal order of emergence is a generically localizable property. Then \cref{Turaev consequence} implies the dissipative counterpart of \cref{typicality conserva emer}:
 \begin{theorem}[\cite{BB19}]
For every $1\le r\le \infty$, a topologically generic $f\in \mathcal N^r$ has its metric emergence of maximal order $2$:
$\mathcal {OE}_\leb= \dim M^2\; .$
\end{theorem}
\begin{remark} Using that high emergence is a generically localizable property, we observe that 
Conjecture \ref{Conjecture} implies in particular Conjecture \ref{positive order emergence conjecture}. \end{remark}

For the second kind of examples, we will use the following notion:
\begin{definition}[\cite{BB19}] \label{topo emergence} The \emph{  topological emergence  $\cE_{top}$} is the covering number of the space of ergodic probabilities  supported by $K$.\end{definition}
We showed: 
\begin{theorem}[\cite{BB19}] \label{exam topo}
 Let $r>1$ and $f\in C^r(M, M)$ which leaves invariant a hyperbolic, locally maximal, compact set $K$. Then $\mathcal O \cE_{top}=\dim\, K $ if one of the two following conditions holds true:
\begin{itemize}
\item the restriction $f|_K$ is conformal.
\item the restriction $f|_K$ is invertible, the differential has determinant 1 and $\dim M=2$.
\end{itemize}
\end{theorem} 

Let us apply this result to the following:
\begin{theorem}[Hofbauer-Keller \cite{HK95}]
There exists uncountably many unimodal maps $f$ of the interval such that for any of its invariant measure  $\mu$, for Lebesgue almost any point $x$, a subsequence of its empirical measures $(\mathsf e_n(x))_n$ converges to $\mu$. 
\end{theorem}
Actually such dynamics turn out to have a hyperbolic, locally maximal, compact set $K$ of box dimension arbitrarily close to $1$. So we deduce:
\begin{corollary}
There exists uncountably many unimodal maps $f$   such that $\cal O\cE_{\leb}(f)=1.$
\end{corollary}
In his PhD Thesis Aminosadat Talebi  is showing the complex counterpart of this result, and moreover that this phenomena is topologically generic in the total bifurcation locus.  

Recently it has been shown the following:
\begin{theorem}[Kiriki-Nakano-Soma \cite{KNS19}]
For every $2\le r<\infty$, if $f\in \Diff^r(M^2)$ displays a saddle periodic point whose local stable and unstable manifold are tangent, then there exists a $C^r$-perturbation $\tilde f$ of $f$ which displays a non-empty domain with the following properties:
\begin{itemize}
\item $\diam f^n(U)\to 0$,
\item For every $x\in U$, the set of cluster values of $(\mathsf e_n(x))_n$ has infinite dimension.
\end{itemize}
In particular $(\mathsf e_n|_U)_n$ diverges $\Leb.$ a.e. and the emergence of $f$ is high.\end{theorem}
In a work in progress with Biebler we are strengthening this theorem to obtain any topology (including $C^\infty$, analytic and polynomial) to obtain local density of emergence with positive order.

\subsection{Toward a theory on dynamics with high emergence}\label{dico entropy and emergence }
Let $X$ be a compact set of finite box dimension and $f$ a continuous map of $X$. 

Here is the emergent counterpart of the variational principle for entropy:
\begin{theorem}[\cite{BB19}]\label{varia emergence}The maximum of the orders of the metric emergences is attained and equal to the order of the topological emergence:
$$\mathcal O\cE_{top}= \max_{\mu\in \cal M(X): f_*\mu=\mu} \mathcal O\cE_{\mu}\; .$$
\end{theorem}
We can reformulate both the metric emergence and the metric entropy in terms of \emph{quantization number}. We recall that the quantization number $\mathcal Q_\mu(\epsilon)$  of a probability measure $\mu$ of a compact metric space $(Y,d)$ is:
$$\mathcal Q_\mu(\epsilon) =\min\left\{N\ge 1 : \exists (x_i)_i \in Y^N\& (a_i)_i\in \R^N\text{ s.t.  } d_W(\mu, \sum_i a_i \delta_{x_i})<\epsilon\right\}\; .$$
The following dictionary is under construction:
\bigskip

\begin{center}

\begin{tabular}{|c|c|}
  \hline
{Topological entropy} (def. \ref{def topological entropy}): &  Topological  emergence (def. \ref{topo emergence}):\\
$\frac 1n\log$(covering $\#$ of $M$ for $d_n$).&covering $\#$ of space of erg. measures.\\
 &\\
  \hline
{ Metric entropy}   of $\nu$ ergodic  \cite[Th.A.7]{BB19}:& { Metric emergence }of $\mu$ invariant\cite[Prop.3.11]{BB19}:\\
$h_\nu=\lim_{n\to \infty} \frac1n \log \mathcal Q^{\mathsf d_n}_\nu$
&$\cE_\mu= \mathcal Q_{\mathsf e_* \mu}$\\
 &\\
  \hline
{Variational Principle}  \ref{varia entropy}:& {Variational Principle} \ref{varia emergence}:\\
$h_{top}=\sup_\nu h_\nu$ & $\mathcal O\cE_{top}= \max_{\mu} \mathcal O\cE_{\mu}$  
  \\
 &\\
  \hline
Margulis-Ruelle inequality  \ref{Ruelle ineq}: &   {Theorem \ref{analogue Ruelle}:} \\
 $h_\nu \le \dim M\cdot \int \log\| Df \| d\nu$  & $\mathcal {OE}_\mu \le \dim \mu$ \\
 &\\
  \hline
Positive entropy conjecture \ref{positive entropy conjecture}:& Conjecture  \ref{positive order emergence conjecture}:\\
Dynamics with positive entropy  are typical &Dynamics with high emergence are typical\\
&\\
  \hline
Ledrappier-Young Entropy formula \ref{Ledrappier-Young entropy formula}. & ? \\
     &     \\
  \hline
  Maximal entropy measure &  \\
   Equidistribution of periodic points &   ? \\
   &\\
  \hline
\end{tabular}
\end{center}
\medskip

To define the counterpart of the uniqueness of the measure with maximal entropy, an idea is to look at a sequence $(\epsilon_n)_n\to 0$ and a sequence $(E_n)_{n\ge 0}$ of maximal $\epsilon_n$-separated finite set in the space of ergodic measures of $f$, and then to look at the convergence of:
\[\hat \mu_n := \frac1{\# E_n} \sum_{\nu\in E_n}\delta_{\nu}\; .\]
In general, this sequence of measures does not need to converge. The same phenomenon happened for the measures equidistributed on the periodic points. For instance, a generic example having a fast growth of the number of periodic points in \cref{theoA} satisfies that the following sequence does not converge:
\[
\left( \frac1{\#\text{Per}_n}  \sum_{x\in \text{Per}_n}\delta_x\right)_n 
\; .\]
\begin{question} 
For the map $\theta\in \R/\Z\mapsto 2\theta\in \R/\Z$, does $(\hat \mu_n )_n$ converge toward a measure $\hat \mu$ for every $(\epsilon_n)_n$ and $(E_n)_{n\ge 0}$ ?
With $\mu := \int_{\cal M_1(X)} \nu d\hat \mu(\nu)$, does $\cal O \cE_\mu=1$ (as  a measure of maximal emergence given by \cref{exam topo,varia emergence} )?
\end{question}
The Ledrappier-Young entropy formula follows from an autosimilarity (which is constant by ergodicity). Roughly speaking, the scaling factor is given by the Lyapunov exponent and the number of branches is given by the entropy. 

Nevertheless, a dynamics with high emergence might have an ergodic decomposition $\mathsf e_*\mu$ which is not self similar. One might need to be more restrictive to obtain autosimilarity. In order to obtain this, it might useful to deal with a local version of the emergence. 

\begin{definition}
The \emph{local order of emergence} at an ergodic probability $\nu$ is 
\[\cal O\cE^{\text{loc}}_\mu (\nu) = \limsup_{\epsilon\to 0}  \frac{\log -\log \mathsf e_*\mu ( B(\nu, \epsilon))}{-\log \epsilon}\; .\]
\end{definition}
\begin{problem}Does the following inequality is correct?
\[\int \cal O\cE^{\text{loc}}_\mu (\nu)d\mathsf e_* \mu(\nu) \le \cal O\cE_\mu\; . 
\]
\end{problem} 
Then an idea would be to consider systems of local order of emergence $\mathsf e_*\mu$-a.e. constant, and then  study the autosimilarity of their ergodic decomposition to deduce the analogous of the entropy formula.

\emph{I am thankful to J.P. Thouvenot for providing me references and to  S. Biebler for his careful reading of this article. I am grateful to the anonymous referee for his advice.}

\bibliographystyle{alpha}
\bibliography{references}

\def\cprime{$'$} \def\cprime{$'$} \def\cprime{$'$} \def\cprime{$'$}
\begin{thebibliography}{MdMvS92}

\bibitem[AM65]{AM65}
M.~Artin and B.~Mazur.
\newblock On periodic points.
\newblock {\em Ann. of Math. (2)}, 81:82--99, 1965.

\bibitem[Arn04]{Ar00}
V.~I. Arnold.
\newblock {\em Arnold's problems}.
\newblock Springer-Verlag, Berlin; PHASIS, Moscow, 2004.
\newblock Translated and revised edition of the 2000 Russian original, With a
  preface by V. Philippov, A. Yakivchik and M. Peters.

\bibitem[{Asa}16]{As16}
M.~{Asaoka}.
\newblock {Fast growth of the number of periodic points in generic families of
  two-dimensional real-analytic area-preserving diffeomorphisms}.
\newblock {\em arXiv}, 1603.08639, March 2016.

\bibitem[AST18]{AST2018}
M.~Asaoka, K.~Shinohara, and D.~Turaev.
\newblock Fast growth of the number of periodic points arising from
  heterodimensional connections.
\newblock {\em arXiv preprint arXiv:1808.07218}, 2018.

\bibitem[Aus18]{A18}
T.~Austin.
\newblock Measure concentration and the weak {P}insker property.
\newblock {\em Publ. Math. Inst. Hautes \'{E}tudes Sci.}, 128:1--119, 2018.

\bibitem[BB19]{BB19}
P.~Berger and J.~Bochi.
\newblock On emergence and complexity of ergodic decompositions.
\newblock {\em arXiv preprint arXiv:1901.03300}, 2019.

\bibitem[BC91]{BC2}
M.~Benedicks and L.~Carleson.
\newblock The dynamics of the {H}\'enon map.
\newblock {\em Ann. Math.}, 133:73--169, 1991.

\bibitem[BCS18]{BCS18}
J.~Buzzi, S.~Crovisier, and O.~Sarig.
\newblock Measures of maximal entropy for surface diffeomorphisms.
\newblock {\em arXiv preprint arXiv:1811.02240}, 2018.

\bibitem[BD99]{BD99}
C.~Bonatti and L.~D{\'{\i}}az.
\newblock Connexions h\'et\'eroclines et g\'en\'ericit\'e d'une infinit\'e de
  puits et de sources.
\newblock {\em Ann. Sci. \'Ecole Norm. Sup. (4)}, 32(1):135--150, 1999.

\bibitem[Ber16]{BE15}
P.~Berger.
\newblock Generic family with robustly infinitely many sinks.
\newblock {\em Invent. Math.}, 205(1):121--172, 2016.

\bibitem[{Ber}17a]{Be16}
P.~{Berger}.
\newblock {Emergence and non-typicality of the finiteness of the attractors in
  many topologies}.
\newblock {\em Proc. of the Steklov Instit.}, {\bf{2}}97(1):1--27, 2017.

\bibitem[Ber17b]{BerPer2018}
P.~Berger.
\newblock Generic family displaying robustly a fast growth of the number of
  periodic points.
\newblock {\em arXiv preprint arXiv:1701.02393}, 2017.

\bibitem[Ber19a]{berhen}
P.~Berger.
\newblock Abundance of non-uniformly hyperbolic {H}\'enon like endomorphisms.
\newblock {\em Ast\'erisque}, 410:53--177, 2019.

\bibitem[Ber19b]{berent}
P.~Berger.
\newblock Properties of the maximal entropy measure and geometry of {H}\'enon
  attractors.
\newblock {\em to appear in JEMS}, 2019.

\bibitem[BGV07]{BGV}
F.~Bolley, A.~Guillin, and C.~Villani.
\newblock Quantitative concentration inequalities for empirical measures on
  non-compact spaces.
\newblock {\em Probab. Theory Related Fields}, 137(3-4):541--593, 2007.

\bibitem[Bie16]{biebler2016persistent}
S.~Biebler.
\newblock Persistent homoclinic tangencies and infinitely many sinks for
  residual sets of automorphisms of low degree in c\^{}$\{$3$\}$.
\newblock {\em arXiv preprint arXiv:1611.02011}, 2016.

\bibitem[BK32]{BK32}
G.~D. Birkhoff and B.~O. Koopman.
\newblock Recent contributions to the ergodic theory.
\newblock {\em Proceedings of the National Academy of Sciences of the United
  States of America}, 18(3):279, 1932.

\bibitem[BLS93]{BLS}
E.~Bedford, M.~Lyubich, and J.~Smillie.
\newblock Polynomial diffeomorphisms of {${\bf C}^2$}. {IV}. {T}he measure of
  maximal entropy and laminar currents.
\newblock {\em Invent. Math.}, 112(1):77--125, 1993.

\bibitem[Bow78]{Bo78}
R.~Bowen.
\newblock {\em On {A}xiom {A} diffeomorphisms}.
\newblock American Mathematical Society, Providence, R.I., 1978.
\newblock Regional Conference Series in Mathematics, No. 35.

\bibitem[BR75]{BR75}
R.~Bowen and D.~Ruelle.
\newblock The ergodic theory of {A}xiom {A} flows.
\newblock {\em Invent. Math.}, 29(3):181--202, 1975.

\bibitem[BT19]{BT19}
P.~Berger and D.~Turaev.
\newblock On herman's positive entropy conjecture.
\newblock {\em Advances in Mathematics}, 2019.

\bibitem[Bur]{Bu17}
D.~Burguet.
\newblock Periodic expansiveness of smooth surface diffeomorphisms and
  applications.
\newblock {\em to appear in JEMS}.

\bibitem[Buz97]{Bu97}
G.~T. Buzzard.
\newblock Infinitely many periodic attractors for holomorphic maps of {$2$}
  variables.
\newblock {\em Ann. of Math. (2)}, 145(2):389--417, 1997.

\bibitem[BY19]{BY19}
P.~Berger and J.-C. Yoccoz.
\newblock Strong regularity.
\newblock {\em Ast\'erisque}, 410:1--13, 2019.

\bibitem[Din71]{Di71}
E.~I. Dinaburg.
\newblock A connection between various entropy characterizations of dynamical
  systems.
\newblock {\em Izv. Akad. Nauk SSSR Ser. Mat.}, 35:324--366, 1971.

\bibitem[DNT16]{DNT16}
T.-C. {Dinh}, V.-A. {Nguyen}, and T.~T. {Truong}.
\newblock {Growth of the number of periodic points for meromorphic maps}.
\newblock {\em ArXiv e-prints}, January 2016.

\bibitem[Dua99]{Duarte99}
P.~Duarte.
\newblock Abundance of elliptic isles at conservative bifurcations.
\newblock {\em Dynam. Stability Systems}, 14(4):339--356, 1999.

\bibitem[FHY83]{FHY83}
A.~Fathi, M.-R. Herman, and J.-C. Yoccoz.
\newblock A proof of {P}esin's stable manifold theorem.
\newblock In {\em Geometric dynamics ({R}io de {J}aneiro, 1981)}, volume 1007
  of {\em Lecture Notes in Math.}, pages 177--215. Springer, Berlin, 1983.

\bibitem[FM89]{Milnor-Friedman}
S.~Friedland and J.~Milnor.
\newblock Dynamical properties of plane polynomial automorphisms.
\newblock {\em Ergodic Theory Dynam. Systems}, 9(1):67--99, 1989.

\bibitem[GHK06]{GHK06}
A.~Gorodetski, B.~Hunt, and V.~Kaloshin.
\newblock Newton interpolation polynomials, discretization method, and certain
  prevalent properties in dynamical systems.
\newblock In {\em International {C}ongress of {M}athematicians. {V}ol. {III}},
  pages 27--55. Eur. Math. Soc., Z{\"u}rich, 2006.

\bibitem[Had98]{Ha98}
J.~Hadamard.
\newblock Sur la forme des lignes g\'{e}od\'{e}siques \`a l'infini et sur les
  g\'{e}od\'{e}siques des surfaces r\'{e}gl\'{e}es du second ordre.
\newblock {\em Bull. Soc. Math. France}, 26:195--216, 1898.

\bibitem[Hed39]{He39}
G.~A. Hedlund.
\newblock The dynamics of geodesic flows.
\newblock {\em Bull. Amer. Math. Soc.}, 45(4):241--260, 1939.

\bibitem[Her83]{He83}
M.-R. Herman.
\newblock {\em Sur les courbes invariantes par les diff\'{e}omorphismes de
  l'anneau. {V}ol. 1}, volume 103 of {\em Ast\'{e}risque}.
\newblock Soci\'{e}t\'{e} Math\'{e}matique de France, Paris, 1983.
\newblock With an appendix by Albert Fathi, With an English summary.

\bibitem[Her98]{He98}
M.~Herman.
\newblock Some open problems in dynamical systems.
\newblock In {\em Proceedings of the {I}nternational {C}ongress of
  {M}athematicians, {V}ol. {II} ({B}erlin, 1998)}, number Extra Vol. II, pages
  797--808, 1998.

\bibitem[HK95]{HK95}
F.~Hofbauer and G.~Keller.
\newblock Quadratic maps with maximal oscillation.
\newblock In {\em Algorithms, fractals, and dynamics ({O}kayama/{K}yoto,
  1992)}, pages 89--94. Plenum, New York, 1995.

\bibitem[Hop39]{Ho39}
E.~Hopf.
\newblock Statistik der geod\"{a}tischen {L}inien in {M}annigfaltigkeiten
  negativer {K}r\"{u}mmung.
\newblock {\em Ber. Verh. S\"{a}chs. Akad. Wiss. Leipzig}, 91:261--304, 1939.

\bibitem[HSY92]{HSY92}
B.~R. Hunt, T.~Sauer, and J.~A. Yorke.
\newblock Prevalence: a translation-invariant ``almost every'' on
  infinite-dimensional spaces.
\newblock {\em Bull. Amer. Math. Soc. (N.S.)}, 27(2):217--238, 1992.

\bibitem[IL99]{IL99}
Y.~Ilyashenko and W.~Li.
\newblock {\em Nonlocal bifurcations}, volume~66 of {\em Mathematical Surveys
  and Monographs}.
\newblock American Mathematical Society, Providence, RI, 1999.

\bibitem[Jak81]{Ja81}
M.~V. Jakobson.
\newblock Absolutely continuous invariant measures for one-parameter families
  of one-dimensional maps.
\newblock {\em Comm. Math. Phys.}, 81(1):39--88, 1981.

\bibitem[Kal99]{K99}
V.~Kaloshin.
\newblock An extension of the {A}rtin-{M}azur theorem.
\newblock {\em Ann. of Math. (2)}, 150(2):729--741, 1999.

\bibitem[Kal00]{K00}
V.~Kaloshin.
\newblock Generic diffeomorphisms with superexponential growth of number of
  periodic orbits.
\newblock {\em Comm. Math. Phys.}, 211(1):253--271, 2000.

\bibitem[Kat80]{Ka80}
A.~Katok.
\newblock Lyapunov exponents, entropy and periodic orbits for diffeomorphisms.
\newblock {\em Inst. Hautes \'{E}tudes Sci. Publ. Math.}, (51):137--173, 1980.

\bibitem[KH]{KH072}
V.~Kaloshin and B.~Hunt.
\newblock Stretched exponential estimates on growth of the number of periodic
  points for prevalent diffeomorphisms {II}.
\newblock {\em http://www.terpconnect.umd.edu/~vkaloshi/papers/Per-pts2.pdf},
  (1):1--177.

\bibitem[KH07]{KH07}
V.~Kaloshin and B.~Hunt.
\newblock Stretched exponential estimates on growth of the number of periodic
  points for prevalent diffeomorphisms. {I}.
\newblock {\em Ann. of Math. (2)}, 165(1):89--170, 2007.

\bibitem[Klo12]{Kl12}
B.~Kloeckner.
\newblock A generalization of {H}ausdorff dimension applied to {H}ilbert cubes
  and {W}asserstein spaces.
\newblock {\em J. Topol. Anal.}, 4(2):203--235, 2012.

\bibitem[KNS19]{KNS19}
S.~Kiriki, Y.~Nakano, and T.~Soma.
\newblock Emergence via non-existence of averages.
\newblock {\em arXiv preprint arXiv:1904.03424}, 2019.

\bibitem[Kol57]{Ko57}
A.~N. Kolmogorov.
\newblock Th\'{e}orie g\'{e}n\'{e}rale des syst\`emes dynamiques et
  m\'{e}canique classique.
\newblock In {\em Proceedings of the {I}nternational {C}ongress of
  {M}athematicians, {A}msterdam, 1954, {V}ol. 1}, pages 315--333, 1957.

\bibitem[KS06]{KS06}
V.~Kaloshin and M.~Saprykina.
\newblock Generic 3-dimensional volume-preserving diffeomorphisms with
  superexponential growth of number of periodic orbits.
\newblock {\em Discrete Contin. Dyn. Syst.}, 15(2):611--640, 2006.

\bibitem[LY85]{LYII85}
F.~Ledrappier and L.-S. Young.
\newblock The metric entropy of diffeomorphisms. {II}. {R}elations between
  entropy, exponents and dimension.
\newblock {\em Ann. of Math. (2)}, 122(3):540--574, 1985.

\bibitem[Lyu83]{Lyu83}
M.~Lyubich.
\newblock Entropy properties of rational endomorphisms of the {R}iemann sphere.
\newblock {\em Ergodic Theory Dynam. Systems}, 3(3):351--385, 1983.

\bibitem[Ma{\~n}82]{Ma82}
R.~Ma{\~n}{\'e}.
\newblock An ergodic closing lemma.
\newblock {\em Ann. of Math. (2)}, 116(3):503--540, 1982.

\bibitem[Ma{\~n}83]{Ma83}
R.~Ma{\~n}{\'e}.
\newblock On the uniqueness of the maximizing measure for rational maps.
\newblock {\em Bol. Soc. Brasil. Mat.}, 14(1):27--43, 1983.

\bibitem[MdMvS92]{MdMvS92}
M.~Martens, W.~de~Melo, and S.~van Strien.
\newblock Julia-{F}atou-{S}ullivan theory for real one-dimensional dynamics.
\newblock {\em Acta Math.}, 168(3-4):273--318, 1992.

\bibitem[Mos62]{Mo62}
J.~Moser.
\newblock On invariant curves of area-preserving mappings of an annulus.
\newblock {\em Nachr. Akad. Wiss. G\"{o}ttingen Math.-Phys. Kl. II},
  1962:1--20, 1962.

\bibitem[New74]{Newhouse}
S.~E. Newhouse.
\newblock Diffeomorphisms with infinitely many sinks.
\newblock {\em Topology}, 12:9--18, 1974.

\bibitem[Pal00]{Pa00}
J.~Palis.
\newblock A global view of dynamics and a conjecture on the denseness of
  finitude of attractors.
\newblock {\em Ast\'erisque}, (261):xiii--xiv, 335--347, 2000.
\newblock G{\'e}om{\'e}trie complexe et syst{\`e}mes dynamiques (Orsay, 1995).

\bibitem[Pal05]{Pa05}
J.~Palis.
\newblock A global perspective for non-conservative dynamics.
\newblock 22(4):485--507, 2005.

\bibitem[Pal08]{Pa08}
J.~Palis.
\newblock Open questions leading to a global perspective in dynamics.
\newblock {\em Nonlinearity}, 21(4):T37--T43, 2008.

\bibitem[Par64]{Pa64}
W.~Parry.
\newblock Intrinsic {M}arkov chains.
\newblock {\em Trans. Amer. Math. Soc.}, 112:55--66, 1964.

\bibitem[Pes76]{Pe76}
J.~B. Pesin.
\newblock Families of invariant manifolds that correspond to nonzero
  characteristic exponents.
\newblock {\em Izv. Akad. Nauk SSSR Ser. Mat.}, 40(6):1332--1379, 1440, 1976.

\bibitem[Poi92]{poincare1892methodes}
H.~Poincar{\'e}.
\newblock {\em Les m{\'e}thodes nouvelles de la m{\'e}canique c{\'e}leste:
  Solutions p{\'e}riodiques. Non-existence des int{\'e}grales uniformes.
  Solutions asymptotiques.-t. 2. M{\'e}thodes de MM. Newcomb, Gyld{\'e}n,
  Lindstedt et Bohlin.-t. 3. Invariants int{\'e}graux. Solutions
  p{\'e}riodiques du deuxi{\`e}me genre. Solutions doublement asymptotiques},
  volume~1.
\newblock Gauthier-Villars, 1892.

\bibitem[Prz82]{Pr82}
F.~Przytycki.
\newblock Examples of conservative diffeomorphisms of the two-dimensional torus
  with coexistence of elliptic and stochastic behaviour.
\newblock {\em Ergodic Theory Dynam. Systems}, 2(3-4):439--463 (1983), 1982.

\bibitem[PS96]{PS95}
C.~Pugh and M.~Shub.
\newblock {\em Stable ergodicity and partial hyperbolicity}, volume 362, pages
  182--187.
\newblock In International Conference on Dynamical Systems (Montevideo, 1995),
  longman, harlow edition, 1996.

\bibitem[PT93]{PT93}
J.~Palis and F.~Takens.
\newblock {\em Hyperbolicity and sensitive chaotic dynamics at homoclinic
  bifurcations}, volume~35 of {\em Cambridge Studies in Advanced Mathematics}.
\newblock Cambridge University Press, Cambridge, 1993.
\newblock Fractal dimensions and infinitely many attractors.

\bibitem[PY09]{PY09}
J.~Palis and J.-C. Yoccoz.
\newblock Non-uniformly hyperbolic horseshoes arising from bifurcations of
  {P}oincar\'e heteroclinic cycles.
\newblock {\em Publ. Math. Inst. Hautes \'Etudes Sci.}, (110):1--217, 2009.

\bibitem[Rue78]{Ru78}
D.~Ruelle.
\newblock An inequality for the entropy of differentiable maps.
\newblock {\em Bol. Soc. Brasil. Mat.}, 9(1):83--87, 1978.

\bibitem[Sin94]{Si94}
Y.~Sina{\u\i}.
\newblock {\em Topics in ergodic theory}, volume~44 of {\em Princeton
  Mathematical Series}.
\newblock Princeton University Press, Princeton, NJ, 1994.

\bibitem[Sma67]{Sm}
S.~Smale.
\newblock Differentiable dynamical systems.
\newblock {\em Bull. Amer. Math. Soc.}, 73:747--817, 1967.

\bibitem[TLY86]{TLY}
L.~Tedeschini-Lalli and J.~A. Yorke.
\newblock How often do simple dynamical processes have infinitely many
  coexisting sinks?
\newblock {\em Comm. Math. Phys.}, 106(4):635--657, 1986.

\bibitem[Tur15]{Tu15}
D.~Turaev.
\newblock Maps close to identity and universal maps in the {N}ewhouse domain.
\newblock {\em Comm. Math. Phys.}, 335(3):1235--1277, 2015.

\bibitem[Ush80]{Us80}
S.~Ushiki.
\newblock Sur les liasons-cols des syst\`emes dynamiques analytiques.
\newblock {\em C. R. Acad. Sci. Paris S\'er. A-B}, 291(7):A447--A449, 1980.

\bibitem[vN43]{vN43}
J.~von Neumann.
\newblock {\em Mathematische {G}rundlagen der {Q}uantenmechanik}.
\newblock Dover Publications, N. Y., 1943.

\bibitem[Whi34]{Wh34}
H.~Whitney.
\newblock Analytic extensions of differentiable functions defined in closed
  sets.
\newblock {\em Trans. Amer. Math. Soc.}, 36(1):63--89, 1934.

\bibitem[Yoc97a]{Y}
J.-C. Yoccoz.
\newblock Jackobson' theorem.
\newblock {\em Manuscript}, 1997.

\bibitem[Yoc97b]{Yolecture1}
J.-C. Yoccoz.
\newblock Quelques exemples de dynamique faiblement hyperbolique.
\newblock {\em
  http://media.college-de-france.fr/media/jean-christophe-yoccoz/UPL3205650719684061902$\_$AN$\_$97$\_$yoccoz.pdf},
  1997.

\bibitem[Yoc16]{Yo16}
J.-C. Yoccoz.
\newblock {H}yperbolicit\'e non-uniforme pour les perturbations
  multidimensionnelles des polyn\^omes quadratiques 1/2.
\newblock {\em
  www.college-de-france.fr/site/jean-christophe-yoccoz/course-2016-02-03-10h00.htm},
  February 2016.

\bibitem[Yoc19]{Y19}
J.~C. Yoccoz.
\newblock A proof of {J}akobson's theorem.
\newblock {\em Ast\'erisque}, 410:15--52, 2019.

\end{thebibliography}

\end{otherlanguage}
\end{document}